\documentclass[11pt]{amsart}
\usepackage[T1]{fontenc}
\usepackage{lmodern,amssymb,mathtools}
\usepackage[margin=28mm]{geometry}
\usepackage{microtype}
\usepackage[hidelinks]{hyperref}
\newtheorem{theorem}{Theorem}[section]
\newtheorem{proposition}[theorem]{Proposition}
\newtheorem{lemma}[theorem]{Lemma}
\newtheorem{corollary}[theorem]{Corollary}
\theoremstyle{definition}
\newtheorem{example}[theorem]{Example}
\theoremstyle{remark}
\newtheorem{remark}[theorem]{Remark}
\numberwithin{equation}{section}
\newcommand{\NN}{\mathbb Z_{\ge0}}
\newcommand{\ZZ}{\mathbb Z}
\newcommand{\QQ}{\mathbb Q}
\newcommand{\CC}{\mathbb C}
\newcommand{\wt}{\operatorname{wt}}

\newcommand{\sg}[1]{\langle #1\rangle}
\newcommand{\qbin}[2]{\begin{bmatrix}#1\\#2\end{bmatrix}_{q}}
\DeclareMathOperator{\Stab}{Stab}
\DeclareMathOperator{\Hilb}{Hilb}
\hypersetup{pdftitle={Positivity of polynomial quotients of q-binomial coefficients},pdfauthor={Shuai Yang and Yueming Zhong}}
\title[Positivity of polynomial $q$-binomial quotients]{Positivity of polynomial quotients of\texorpdfstring{\\}{} $q$-binomial coefficients}
\author[S.~Yang]{Shuai Yang}
\address{Key Laboratory of Engineering Modeling and Statistical Computation of Hainan Province; School of Mathematics and Statistics, Hainan University, Haikou 570228, P.R. China}
\email{yangs0901@126.com}
\author[Y.~Zhong]{Yueming Zhong\textsuperscript{*}}
\address{Key Laboratory of Engineering Modeling and Statistical Computation of Hainan Province; School of Mathematics and Statistics, Hainan University, Haikou 570228, P.R. China}
\email{zhongyueming107@gmail.com}
\thanks{\textsuperscript{*}Corresponding author: Yueming Zhong. Email: \href{mailto:zhongyueming107@gmail.com}{\nolinkurl{zhongyueming107@gmail.com}}.}
\date{27 November 2026}
\begin{document}
\begin{abstract}
We prove the conjecture of Gatzweiler and Krattenthaler that a polynomial
quotient of two Gaussian binomial coefficients with the same upper index
has nonnegative integer coefficients. For $0\le\ell<k\le n/2$, such a
quotient $\qbin nk/\qbin n\ell$ is the Hilbert series of an Artinian
weighted complete intersection, with weights and defining degrees
prescribed by its product formula. The underlying combinatorial result
concerns two disjoint integer intervals of equal length: divisibility
inequalities imply the semigroup inequalities that guarantee the existence
of the required regular sequence. The proof uses Kneser's addition
theorem and an elementary lemma locating a multiple of the largest
generator. In the final case, we change the cyclic modulus from the
largest generator to the smallest. We also give an explicit specialization
over the rational numbers and an algorithm that constructs finite monomial
sets, closed under taking divisors, whose weight enumerators are these
quotients. A four-variable example illustrates the complete intersection
realization.
\end{abstract}
\maketitle

\section{Introduction}
Let $q$ be an indeterminate. For integers $0\le j\le n$, the Gaussian
binomial coefficient is
\[
 \qbin nj=\prod_{i=1}^{j}\frac{1-q^{n-j+i}}{1-q^i},
\]
where an empty product is $1$. Gatzweiler and Krattenthaler
\cite[Conjecture~1]{GK} conjectured that a polynomial quotient
\[
 Q_{n,k,\ell}(q)=\frac{\qbin nk}{\qbin n\ell}
\]
always has nonnegative coefficients. Polynomiality alone does not
imply positivity for arbitrary products of this type; for example,
\[
 \frac{(1-q)(1-q^6)}{(1-q^2)(1-q^3)}=1-q+q^2.
\]
The conjecture concerns the additional structure imposed by the
common upper index.

In the normalized range $0\le\ell<k\le n/2$, Gatzweiler and
Krattenthaler review the previously known cases $\ell\le1$ and
$k-\ell\le2$, and prove the further cases $\ell=2$ and $k-\ell=3$.
Gatzweiler, Levic\'an-Santib\'a\~nez and Yoshida developed a
complete-intersection approach in~\cite{GLY}. Their application uses a
semigroup criterion for the existence of homogeneous regular sequences
with prescribed weights and degrees. An earlier criterion from which
this one follows appears in Bender and Spaenlehauer~\cite[Theorem~3.8]{BenderSpaenlehauer2024};
we explain the reduction in Section~\ref{sec:algebra}. As an application,
\cite[Corollary~3.10]{GLY} proves the conjecture under the condition
\[
 n>2k\left\lfloor\frac{k-1}{2}\right\rfloor+k-\ell.
\]
They also suggest using the interval structure to obtain further cases
from their criterion; see~\cite[Remark~3.12]{GLY}. We establish the
required semigroup inequalities for these intervals throughout the
normalized range. This gives the following result.

\begin{theorem}\label{thm:main}
Let $n,k,\ell$ be integers with $0\le k,\ell\le n$.
If
\[
 Q_{n,k,\ell}(q)=\frac{\qbin nk}{\qbin n\ell}
\]
is a polynomial, then it belongs to $\NN[q]$.
More precisely, if $0\le\ell<k\le n/2$, every such polynomial
is the Hilbert series of an Artinian weighted complete intersection
with weights $\ell+1,\ldots,k$ and defining degrees
$n-k+1,\ldots,n-\ell$.
\end{theorem}

The combinatorial statement underlying Theorem~\ref{thm:main} is
independent of $q$. For integers $u\le v$, write
$[u,v]=\{u,u+1,\ldots,v\}$. For a nonempty finite set $S$ of positive integers, write
\[
 \sg S=\left\{\sum_{s\in S}c_s s:c_s\in\NN\right\}.
\]
Thus $0\in\sg S$.

\begin{theorem}\label{thm:interval-semigroup}
Let $W=[L+1,K]$ with integers $0\le L<K$, put $r=K-L$, and let
$I=[b,b+r-1]$, where $b$ is an integer with $b>K$. If
\begin{equation}\label{eq:interval-D}
 |I\cap d\ZZ|\ge|W\cap d\ZZ|\qquad(d=1,2,\ldots),
\end{equation}
then every nonempty subset $S\subseteq W$ satisfies
\[
 |I\cap\sg S|\ge|S|.
\]
\end{theorem}

The proof works with sumsets in a cyclic group. A consequence of
Kneser's theorem gives either the required semigroup count or a lower
bound for the size of an iterated sumset. With the largest generator as
modulus, translation by the smallest generator settles one range. In the
remaining range, we
locate a multiple of the largest generator in the target interval. Its
position either places a translate of $S$ inside that interval or forces
its multiplier to be large enough for the final count, taken modulo the
smallest generator. The interval structure of $W$ is used to control
this position and to bound the span of $S$.

We also record constructive consequences of the complete-intersection
realization. A determinant argument specifies integer coefficients for
a regular sequence over $\QQ$. A Gr\"obner basis computation then gives
a finite set of monomials, closed under taking divisors, whose weight
enumerator is the quotient.

Section~\ref{sec:algebra} gives the algebraic reduction.
Section~\ref{sec:cyclic} proves the cyclic estimates, and
Section~\ref{sec:interval-proof} proves the two main theorems.
Section~\ref{sec:models} gives the constructive consequence and a
four-variable example. Appendix~\ref{app:criterion} includes a proof
of the sufficient algebraic criterion used here, and
Appendix~\ref{app:effective} gives the explicit specialization.

\section{Algebraic reduction}\label{sec:algebra}
Let $F$ be a field, and let $R=F[x_1,\ldots,x_r]$ be a polynomial ring
with positive weights
$\deg x_i=w_i$. An \emph{Artinian weighted complete intersection}
is a finite-dimensional quotient $R/(f_1,\ldots,f_r)$, where
$f_1,\ldots,f_r$ is a homogeneous regular sequence.
For a graded $F$-algebra $A=\bigoplus_{d\ge0}A_d$ with
finite-dimensional homogeneous components, write
\[
 \Hilb(A;q)=\sum_{d\ge0}(\dim_F A_d)q^d.
\]
For a monomial $x^\alpha=\prod_i x_i^{\alpha_i}$, where
$\alpha\in\NN^r$, its weight is
$\wt(x^\alpha)=\sum_i w_i\alpha_i$.

We use the following sufficient form of the semigroup criterion,
as stated and applied to positivity by
Gatzweiler--Levic\'an-Santib\'a\~nez--Yoshida.
We state it for distinct weights and degrees, which is the case
needed in this paper.

\begin{proposition}[{\cite[Propositions~3.3 and~3.6]{GLY}}]\label{prop:ci}
Let $W,V$ be finite sets of positive integers with $|W|=|V|$.
Suppose that
\begin{equation}\label{eq:H}
 |V\cap\sg S|\ge |S|\qquad(\varnothing\ne S\subseteq W).
\end{equation}
Then there exists a homogeneous regular sequence $(f_v)_{v\in V}$
in $\CC[x_w:w\in W]$, with $\deg x_w=w$ and $\deg f_v=v$.
Its Artinian quotient has Hilbert series
\begin{equation}\label{eq:hilbert-ci}
 \Hilb\bigl(\CC[x_w:w\in W]/(f_v:v\in V);q\bigr)
 =\frac{\prod_{v\in V}(1-q^v)}{\prod_{w\in W}(1-q^w)}\in\NN[q].
\end{equation}
\end{proposition}

The regular-sequence criterion in~\cite[Proposition~3.3]{GLY} is
an equivalence and allows repeated weights and degrees. Its sufficient
direction also follows from~\cite[Theorem~3.8, pp.~178--179]{BenderSpaenlehauer2024}.
To see this when $m=|W|=|V|\ge1$, put $R=\CC[x_w:w\in W]$ and adjoin
a variable $y$ of integer weight $h>\max V$. Apply that theorem with its
parameters $n=r=m$, so its lower bound for a subset $J$ of the variables
is $|J|-1$. A monomial of a degree in $V$ cannot involve $y$.
For $\varnothing\ne S\subseteq W$, taking
$J=\{x_w:w\in S\}\cup\{y\}$ gives the condition
$|V\cap\sg S|\ge|S|$; taking $J=\{x_w:w\in S\}$ gives the weaker
inequality $|V\cap\sg S|\ge|S|-1$. The cases $J=\varnothing$ and
$J=\{y\}$ are automatic. Thus \eqref{eq:H} supplies
all the required inequalities. The resulting
regular sequence in $R[y]$ consists of forms in $R$, where it is also
regular by faithful flatness. The case $m=0$ is immediate.
Appendix~\ref{app:criterion} gives a direct proof of the sufficient
direction. The Hilbert-series identity also follows from
the multiplication exact sequences; see~\cite[Corollary~3.3]{Stanley1978}.

In the normalized range put
\begin{equation}\label{eq:intervals}
 r=k-\ell,\qquad b=n-k+1,\qquad
 W=[\ell+1,k],\qquad V=[b,b+r-1].
\end{equation}
Then $b>k$ and
\begin{equation}\label{eq:product}
 Q_{n,k,\ell}(q)=
 \frac{\prod_{v\in V}(1-q^v)}{\prod_{w\in W}(1-q^w)}.
\end{equation}
Cyclotomic factorization makes polynomiality equivalent to
\begin{equation}\label{eq:D}
 |V\cap d\ZZ|\ge|W\cap d\ZZ|\qquad(d=1,2,\ldots).
\end{equation}
Indeed, the difference of the two counts is the exponent of the
$d$th cyclotomic polynomial. At $d=1$ it is zero. Our task is to
deduce \eqref{eq:H} from \eqref{eq:D} using the two intervals.

\section{Cyclic interval estimates}\label{sec:cyclic}
For a finite nonempty subset $A$ of an abelian group $G$ and an integer
$p\ge1$, let $pA$ denote the sum of $p$ copies of $A$.
For a finite nonempty subset $C\subseteq G$, write
\[\Stab(C)=\{g\in G:C+g=C\}.\]
Kneser's original work~\cite[pp.~463--464]{Kneser1953}
establishes an asymptotic-density theorem for subsets of $\ZZ$.
We use the following iterated consequence of the finite-set
addition theorem in~\cite[Theorem~1]{DeVos2014}:
\begin{equation}\label{eq:kneser}
 |pA|\ge p|A+H|-(p-1)|H|,
 \qquad H=\Stab(pA),\quad p\ge2.
\end{equation}
To obtain \eqref{eq:kneser}, let $\bar A$ be the image of $A$ in $G/H$.
A period of $j\bar A$ is also a period of $p\bar A$, so
\[
 \Stab(j\bar A)\subseteq\Stab(p\bar A)=\{0\}
 \qquad(1\le j\le p).
\]
Iterating the two-set inequality in $G/H$ gives
$|p\bar A|\ge p|\bar A|-(p-1)$.
Multiplication by $|H|$ yields \eqref{eq:kneser}.

\begin{lemma}[Lifting a residue]\label{lem:lifting}
Let $S$ be a nonempty finite set of positive integers, let $m\in S$,
and let $J=[b,b+m-1]$ with an integer $b\ge1$.
If $x\in\sg S$ and $0\le x\le b+m-1$, the unique representative
$y\in J$ of $x$ modulo $m$ also belongs to $\sg S$.
\end{lemma}
\begin{proof}
The difference $x-y$ is a multiple of $m$ and is at most $m-1$.
Consequently $y=x+jm$ for some integer $j\ge0$.
\end{proof}

\begin{lemma}[Growth or a sufficient count]\label{lem:growth}
Let $S$ be a nonempty finite set of positive integers, let $m\in S$,
and assume that reduction modulo $m$ is injective on $S$.
Write $t=|S|$ and $A=S\bmod m\subseteq\ZZ/m\ZZ$.
Let $I=[b,b+r-1]$, where $b,r$ are integers with $b\ge1$ and $1\le r<m$, and suppose
\begin{equation}\label{eq:Dlocal}
 |I\cap d\ZZ|\ge|S\cap d\ZZ|\qquad(d\mid m).
\end{equation}
Put $J=[b,b+m-1]$. If $p\ge2$ is an integer and every residue of $C=pA$
has a representative in $J\cap\sg S$, then either
\[
 |I\cap\sg S|\ge t\qquad\text{or}\qquad |C|\ge p(t-1)+1.
\]
\end{lemma}
\begin{proof}
We have $|A|=t$ and $0\in A$. Set
\[
 H=\Stab(C),\quad h=|H|,\quad
 \mu=|A+H|/h,\quad e=\mu h-t.
\]
Thus $e\ge0$, and Kneser's theorem gives
\begin{equation}\label{eq:coset-growth}
 |C|/h\ge1+p(\mu-1).
\end{equation}
Let $d=m/h$. The subgroup $H$ consists of the multiples of $d$
modulo $m$, and injectivity implies
$|S\cap d\ZZ|=|A\cap H|$.

If $\mu=1$, then $A\subseteq H$ since $0\in A$, and $C\subseteq H$.
As $C$ is a nonempty union of $H$-cosets, $C=H$.
Condition \eqref{eq:Dlocal} supplies at least $t$ multiples of $d$
in $I$. Since $C=H$ and $I\subseteq J$, the hypothesis on the
representatives of $C$ shows that all these multiples belong to $\sg S$.
Hence the first alternative holds.

Suppose $\mu\ge2$. Let $T$ be the set of residues represented by the
omitted tail
$J\setminus I$, and write
\[
 e_0=h-|A\cap H|\le e.
\]
Here $e$ counts the points missing from $A$ in the $\mu$ cosets
forming $A+H$, and $e_0$ counts those missing from the single coset $H$.
There are $h$ representatives of $H$ in $J$.
The inequality at $d=m/h$ shows that at most $e_0$ of them lie in
$J\setminus I$. Since this tail is a consecutive integer interval,
its counts in distinct residue classes modulo $d$ differ by at most
one. Thus every other $H$-coset loses at most $e_0+1\le e+1$
representatives when $J$ is restricted to $I$.

If $(p-1)h\ge p(e+1)$, then $p(h-e-1)\ge h>0$, so $h-e-1\ge1$.
Since $0\in C$ and $C+H=C$, the set $C$ contains $H$. By
\eqref{eq:coset-growth}, it contains at least
$p(\mu-1)$ other $H$-cosets. The representatives of $C$ that remain
in $I$ therefore number at least
\[
 (h-e)+p(\mu-1)(h-e-1)
 =t+(\mu-1)((p-1)h-p(e+1))\ge t.
\]
Otherwise, since all quantities are integers,
$pe-(p-1)h\ge-p+1$. Kneser's inequality now gives
\[
 |C|\ge p(t+e)-(p-1)h\ge p(t-1)+1,
\]
which is the second alternative.
\end{proof}

\subsection{Translation by the smallest generator}
\begin{lemma}\label{lem:translation}
Let $S$ be a nonempty finite set of positive integers, and put
$a=\min S$, $M=\max S$, $t=|S|$. Let $b,r$ be positive integers with $b>M$, put
$I=[b,b+r-1]$,
and suppose
\[
 |I\cap d\ZZ|\ge |S\cap d\ZZ|\qquad(d\mid M).
\]
If
\begin{equation}\label{eq:translation-range}
 M-r\le a\le r,
\end{equation}
then $|I\cap\sg S|\ge t$.
\end{lemma}
\begin{proof}
If $r\ge M$, lift the $t$ distinct residues of $S$ modulo $M$
to the first $M$ integers of $I$ using Lemma~\ref{lem:lifting}.
Assume $r<M$ and put
\[
 \delta=M-r,\qquad J=[b,b+M-1],\qquad A=S\bmod M.
\]
The residues $T$ of $J\setminus I$ form a cyclic interval of
length $\delta$. Since \eqref{eq:translation-range} says
$\delta\le a\le M-\delta$, we have
\begin{equation}\label{eq:disjoint}
 T\cap(T+a)=\varnothing.
\end{equation}
Indeed, two members of an interval of length $\delta$ cannot differ,
modulo $M$, by an integer in $[\delta,M-\delta]$.
Also $t\ge2$, because a singleton would have $a=M>r$.

First suppose $b\ge2M$. Take $p=\lfloor b/M\rfloor\ge2$ and
$C=pA$. Represent the zero residue of $A$ by zero and every other
residue by the corresponding member of $S$. Each residue of $C$
then has a semigroup representative at most $p(M-1)$.
Those of $a+C$ have representatives at most
\[
 a+p(M-1)\le pM+M-1\le b+M-1.
\]
Thus every residue of either sumset has a representative in $J\cap\sg S$.
If $I$ contained at most $t-1$ semigroup elements,
Lemma~\ref{lem:growth} would imply $|C|\ge2t-1$.
Since $|C\setminus T|\le t-1$, it follows that $|C\cap T|\ge t$.
By \eqref{eq:disjoint}, translation by $a$ sends these at least $t$
residues outside $T$. Their representatives in $J$ belong to $\sg S$
because the translated residues lie in $a+C$. They therefore give $t$
semigroup elements in $I$, a contradiction.

Now suppose $M<b<2M$. Write $b=M+s$, where $1\le s<M$, and
take $C=2A$. Each residue of $C$ has a semigroup representative at most
$2(M-1)$, so its representative in $J$ also lies in $\sg S$.
The condition at $d=M$ forces $2M\in I$;
there is no other multiple of $M$ that could belong to this interval.
Consequently $s\ge\delta+1$. Identifying residues with their
representatives in $[0,M-1]$, we have
\[
 T=[s-\delta,s-1]\subseteq[1,M-1].
\]
Suppose that $|I\cap\sg S|\le t-1$. Since $A\subseteq C$, $|A|=t$,
and every residue of $C$ lifts to $J\cap\sg S$, we have
$A\cap T\ne\varnothing$. Every member of $A\cap T$
is a positive member of $S$, so $a\le s-1$.
The residues of $a+C$ therefore also have representatives in
$J\cap\sg S$, since they have semigroup representatives at most
\[
 a+2(M-1)\le2M+s-3<b+M-1.
\]
Under the same assumption, Lemma~\ref{lem:growth} gives
$|C|\ge2t-1$, so $|C\cap T|\ge t$. By \eqref{eq:disjoint}, translation
by $a$ sends these points into $(a+C)\setminus T$. Their representatives
give at least $t$ semigroup elements in $I$, again a contradiction.
\end{proof}

\section{The interval theorem and positivity}\label{sec:interval-proof}
\subsection{The position of a multiple}
The next lemma uses every integer between the smallest and largest generators.

\begin{lemma}\label{lem:position}
Let $r<a<M$ and $b>M$ be positive integers. Suppose that
$I=[b,b+r-1]$ contains a multiple of every integer in $[a,M]$.
Let $pM$ be its unique multiple of $M$, and put
\[
 z=pM-b,\qquad w=M-a,\qquad E=M-r.
\]
Then $p\ge2$, $0\le z\le r-1$, and
\begin{equation}\label{eq:position}
 z\ge w\qquad\text{or}\qquad p\ge E.
\end{equation}
\end{lemma}
\begin{proof}
The existence of $pM$ is a hypothesis, and $r<M$ makes it unique.
Since $b>M$, we have $p\ge2$.
Suppose $z<w$, and set
\[
 j=\left\lfloor\frac zp\right\rfloor+1,\qquad d=M-j.
\]
Then $1\le j\le z+1\le w$, so $d\in[a,M-1]$.
Also $pj>z$, which gives
\[
 pd=pM-pj<pM-z=b.
\]
Since a multiple of $d$ belongs to $I$, its multiplier is at least
$p+1$. Necessarily
\[
 (p+1)(M-j)\le b+r-1=pM-z+r-1,
\]
or
\begin{equation}\label{eq:gap}
 E+z+1\le(p+1)j.
\end{equation}
Write $z=p(j-1)+v$, where $0\le v\le p-1$.
Inequality \eqref{eq:gap} becomes $j\ge E-p+v+1$, and hence
\[
 z=p(j-1)+v\ge p(E-p)+(p+1)v\ge p(E-p).
\]
If $p\le E-1$, then
\[
 p(E-p)-(E-1)=(p-1)(E-p-1)\ge0.
\]
But $r<a$ gives $w=M-a\le E-1$, so this would imply
$z\ge E-1\ge w$, a contradiction. Thus $p\ge E$.
\end{proof}

\subsection{Proofs of the main theorems}

\begin{proof}[Proof of Theorem~\ref{thm:interval-semigroup}]
Fix $S$ and write $t=|S|$, $a=\min S$, and $M=\max S$.
The inequalities for $W$ imply those for $S$.
For $t=1$, use \eqref{eq:interval-D} at $d=M$.
Assume $t\ge2$.

If $M\le r$, the first $M$ integers of $I$ contain the $t$ distinct
residues of $S$ modulo $M$. Their representatives are in $\sg S$
by Lemma~\ref{lem:lifting}. Hence assume $M>r$.
Since $M\le K$ and $a\ge L+1$,
\[
 M-r\le L<a.
\]
If $a\le r$, Lemma~\ref{lem:translation} applies.

It remains to consider $a>r$. Let $w=M-a$ and $E=M-r$.
Since $W$ is an interval of length $r$,
\begin{equation}\label{eq:span-short}
 t-1\le w\le r-1<a,\qquad M<2a.
\end{equation}
Every $d\in[a,M]$ belongs to $W$, so \eqref{eq:interval-D}
guarantees that $I$ contains a multiple of $d$.
Let $pM\in I$ and $z=pM-b$ be as in Lemma~\ref{lem:position}.

If $z\ge w$, the $t$ distinct integers
\[
 (p-1)M+S\subseteq[pM-w,pM]\subseteq I
\]
belong to $\sg S$, giving the required count.

Otherwise Lemma~\ref{lem:position} gives $p\ge E$.
Now change the modulus to $a$. By \eqref{eq:span-short}, reduction
modulo $a$ is injective on $S$. Put
\[
 A=S\bmod a,\qquad C=pA,\qquad J=[b,b+a-1].
\]
Each residue of $C$ has a representative that is a sum of
$p$ members of $S$, and is therefore at most $pM$.
Since $pM\in I$ and $r<a$, this is at most $b+a-1$.
Lemma~\ref{lem:lifting} shows that every residue of $C$ is represented
in $J\cap\sg S$. The hypotheses of Lemma~\ref{lem:growth} with
modulus $a$ are now satisfied: in particular, for every $d\mid a$,
\[
 |I\cap d\ZZ|\ge |W\cap d\ZZ|\ge |S\cap d\ZZ|.
\]

If the desired count failed, that lemma would give
$|C|\ge p(t-1)+1$. Removing the $a-r$ residues of $J\setminus I$
would leave at least
\begin{align*}
 |C|-(a-r)
 &\ge p(t-1)+1-(E-w)\\
 &\ge p(t-1)+t-E\\
 &\ge t
\end{align*}
represented points in $I$. The second inequality uses $w\ge t-1$;
the last follows from $p(t-1)-E\ge p-E\ge0$, since $t\ge2$ and $p\ge E$.
This contradiction completes
the proof.
\end{proof}

\begin{remark}
Two features of the interval $W$ enter the proof.
Its length gives $M<2a$ in the last case, ensuring injective reduction
modulo $a$. The presence of every intermediate integer ensures that $d=M-j$ in
Lemma~\ref{lem:position} belongs to $W$, even when $d\notin S$.
The divisibility conditions for $S$ alone would not justify this step.
\end{remark}

\begin{proof}[Proof of Theorem~\ref{thm:main}]
First suppose $0\le\ell<k\le n/2$ and use the intervals
\eqref{eq:intervals}. Polynomiality of \eqref{eq:product} gives
\eqref{eq:D}. Apply Theorem~\ref{thm:interval-semigroup} with
$L=\ell$, $K=k$, and $I=V$. Thus \eqref{eq:H} holds for every
nonempty $S\subseteq W$.
Proposition~\ref{prop:ci} supplies a homogeneous regular sequence
with weights $W$ and degrees $V$. Its Artinian quotient has
Hilbert series $Q_{n,k,\ell}(q)$, proving nonnegativity and the
stronger realization statement.

For general $0\le k,\ell\le n$, use
$\qbin nj=\qbin n{n-j}$ to replace each lower index by its minimum
with its complement. We may therefore assume $k,\ell\le n/2$.
If $k=\ell$, the quotient is $1$.
If $k<\ell$, then
\[
 \deg\qbin nk=k(n-k)<\ell(n-\ell)=\deg\qbin n\ell,
\]
so the nonzero quotient cannot be a polynomial.
The remaining case is the normalized one already proved.
\end{proof}

\section{Standard monomials}\label{sec:models}
Appendix~\ref{app:effective} gives a deterministic choice of a regular
sequence. Its standard monomials provide the following model for the quotient.

\begin{corollary}\label{cor:algorithm}
For every normalized
polynomial quotient
\[
 Q_{n,k,\ell}(q)=\frac{\qbin nk}{\qbin n\ell},
 \qquad 0\le\ell<k\le n/2,
\]
there is a finite deterministic algorithm, depending only on
$(n,k,\ell)$, producing a finite set $\mathcal B_{n,k,\ell}$ of monomials in the variables $x_{\ell+1},\ldots,x_k$, with
$\wt(x_w)=w$, such that
\[
 Q_{n,k,\ell}(q)=\sum_{m\in\mathcal B_{n,k,\ell}}q^{\wt(m)}.
\]
The set is closed under taking monomial divisors, and every monomial
in the set has weight at most $(k-\ell)(n-k-\ell)$.
\end{corollary}
\begin{proof}
Set $W=\{\ell+1,\ldots,k\}$ and
$V=\{n-k+1,\ldots,n-\ell\}$.
Check polynomiality by the cyclotomic inequalities for $1\le d\le k$;
for $d>k$, the denominator count $|W\cap d\ZZ|$ is zero, so the
remaining inequalities are automatic.
Theorem~\ref{thm:main} gives a regular sequence with weights $W$ and
degrees $V$. Proposition~\ref{prop:specialization} specifies
integer-coefficient forms over $\QQ$, using a fixed enumeration of
their monomials. Let $J$ be the ideal they generate in
$R=\QQ[x_w:w\in W]$. Compute its reduced Gr\"obner basis in
lexicographic order $x_{\ell+1}>\cdots>x_k$. This computation terminates by
\cite[Chapter~2, \S7, Theorems~2 and~5]{CoxLittleOShea2015}.
Since $R/J$ is finite dimensional, for every $w\in W$ the basis has
a leading monomial $x_w^{N_w}$ with $N_w\ge1$;
see~\cite[Chapter~5, \S3, Theorem~6]{CoxLittleOShea2015}.
Enumerate the monomials $\prod_{w\in W}x_w^{a_w}$ with
$0\le a_w<N_w$, and retain those divisible by none of the leading
monomials. This finite procedure produces exactly the standard monomials.

The computation preserves weighted homogeneity: the two multiples
forming an $S$-polynomial have the weight of the least common multiple
of their leading monomials, and each division step subtracts a homogeneous multiple of
the same weight as the term being reduced. Thus the reduced basis is
weighted homogeneous, and the normal form of a homogeneous polynomial
is zero or homogeneous of the same weight.
By~\cite[Chapter~5, \S3, Proposition~4]{CoxLittleOShea2015}, the standard
monomials form a basis of $R/J$; they therefore form a homogeneous basis.
Their weight enumerator is the Hilbert series of $R/J$, which proves
the required identity. Their divisor-closed property follows because
the complement of a monomial ideal is closed under divisors.
The degree bound is
\[
 D=\sum_{v\in V}v-\sum_{w\in W}w
  =(k-\ell)(n-k-\ell).
\]
\end{proof}

\begin{example}\label{ex:four}
Consider $(n,k,\ell)=(18,7,3)$. In variables $x,y,z,t$ of weights
$4,5,6,7$, respectively, take
\begin{equation}\label{eq:four-sequence}
 z^2+x^3,\qquad zt+x^2y,\qquad t^2+x^2z,\qquad y^3.
\end{equation}
These forms have degrees $12,13,14,15$ and are a regular sequence
over every field. Indeed, at a common zero over an algebraic closure,
$y=0$ and $zt=0$. If $z=0$, the first and third equations force
$x=t=0$. If $t=0$, the third equation gives $x^2z=0$, and the first
then forces $x=z=0$. The forms are consequently a homogeneous system
of parameters in the polynomial ring over the algebraic closure,
and hence a regular sequence. Regularity descends to the original
field by faithful flatness.

For lexicographic order $t>z>y>x$, the standard monomials split into
the following disjoint blocks:
\begin{align*}
 \mathcal B_0={}&\{x^a:0\le a\le8\}
 \cup\{x^ay:0\le a\le5\}
 \cup\{x^ay^2:0\le a\le4\},\\
 \mathcal B_z={}&\{zx^ay^b:0\le a\le4,\ 0\le b\le1\}
 \cup\{zx^ay^2:0\le a\le1\},\\
 \mathcal B_t={}&\{tx^ay^b:0\le a\le1,\ 0\le b\le2\}
 \cup\{tx^2\}.
\end{align*}
Their sizes are $20$, $12$, and $7$. The supplementary material
gives a monic Gr\"obner basis over $\ZZ$ and explicit ideal-membership
identities. Thus these monomials form a $\ZZ$-basis of
\[
 \ZZ[x,y,z,t]/(z^2+x^3,zt+x^2y,t^2+x^2z,y^3).
\]
With
$[m]_{q^4}=1+q^4+\cdots+q^{4(m-1)}$, their weight enumerator yields
\begin{align*}
 Q_{18,7,3}(q)
 ={}&[9]_{q^4}+q^5[6]_{q^4}+q^{10}[5]_{q^4}\\
 &+q^6\bigl((1+q^5)[5]_{q^4}+q^{10}[2]_{q^4}\bigr)\\
 &+q^7\bigl((1+q^5+q^{10})[2]_{q^4}+q^8\bigr).
\end{align*}
This is a positive decomposition of degree $32$ with total coefficient
sum $39$.
\end{example}

The construction in Appendix~\ref{app:effective} produces regular
sequences over $\QQ$. In Example~\ref{ex:four}, the displayed equations
are regular over every field, and the monic certificate proves that
the integral quotient is free. The general construction does not give
a direct formula for the standard monomials in terms of $(n,k,\ell)$.

\appendix
\section{The sufficient semigroup criterion}\label{app:criterion}
We give a coefficient-space proof of the sufficient direction of the
criterion in~\cite[Proposition~3.3]{GLY}.

\begin{proof}[Proof of Proposition~\ref{prop:ci}]
Put $R=\CC[x_w:w\in W]$, $r=|W|=|V|$, and let $R_v$ be its
weighted degree-$v$ part. The case $r=0$ is immediate. Condition
\eqref{eq:H} with $S=W$ implies $R_v\ne0$ for every $v\in V$.
Consider the affine coefficient space
\[
 \mathcal P=\prod_{v\in V}R_v,\qquad N=\dim\mathcal P.
\]
Fix $\varnothing\ne S\subseteq W$, put
$a_S=|V\cap\sg S|$, and choose $s_0\in S$. Define the torus slice
\[
 Y_S=\{x:x_w=0\ (w\notin S),\ x_w\ne0\ (w\in S),\ x_{s_0}=1\}
       \cong(\CC^*)^{|S|-1}.
\]
At every $x\in Y_S$, evaluation $R_v\to\CC$ is the zero functional
if $v\notin\sg S$, and is nonzero otherwise: in the latter case a
degree-$v$ monomial supported on $S$ does not vanish at $x$.
Consequently the incidence variety
\[
 \mathcal Z_S=\{(f,x)\in\mathcal P\times Y_S:
                       f_v(x)=0\text{ for every }v\in V\}
\]
is a vector bundle of rank $N-a_S$ over $Y_S$. Indeed, one may
choose one such monomial for each $v\in V\cap\sg S$ and solve for its
coefficient, since its value is everywhere invertible on $Y_S$;
the equations use separate coefficient blocks. Hence
\[
 \dim\mathcal Z_S=N+|S|-1-a_S\le N-1.
\]
The closure $B_S$ of its image in $\mathcal P$ is therefore a proper
closed subset. Every common zero with support $S$ can be moved into
$Y_S$ by the weighted action $x_w\mapsto\lambda^w x_w$: choose
$\lambda\in\CC^*$ with $\lambda^{s_0}x_{s_0}=1$.
Homogeneity ensures that all the forms still vanish after this rescaling.
Thus a coefficient
tuple outside $B_S$ has no common zero with support $S$.

There are only finitely many nonempty $S$, and the affine space
$\mathcal P$ is irreducible. Choose
$f\in\mathcal P\setminus\bigcup_{\varnothing\ne S\subseteq W}B_S$.
Its only common zero is the origin. By the Nullstellensatz,
$R/(f_v:v\in V)$ is finite dimensional. The $r$ homogeneous
elements form a system of parameters after localization at the
homogeneous maximal ideal $\mathfrak m=(x_w:w\in W)$. Since
$R_{\mathfrak m}$ is Cohen--Macaulay, they form a regular sequence
there; see~\cite[Tags 00ND and 02JN]{Stacks}. Fix an ordering
$v_1,\ldots,v_r$ of $V$. Each successive homogeneous quotient
$R/(f_{v_1},\ldots,f_{v_j})$, for $0\le j\le r$, embeds in its
localization at $\mathfrak m$: multiplication by a
polynomial with nonzero constant term cannot annihilate the lowest
nonzero homogeneous component of an element. Regularity therefore
descends from $R_{\mathfrak m}$ to $R$.
The successive exact sequences for multiplication by these elements
give
\[
 \operatorname{Hilb}(R/(f_v:v\in V);q)
 =\frac{\prod_{v\in V}(1-q^v)}{\prod_{w\in W}(1-q^w)}.
\]
Each coefficient on the left is a dimension, hence is a nonnegative
integer, and only finitely many are nonzero.
\end{proof}
\section{Effective specialization}\label{app:effective}
\label{sec:specialization}

The following specialization argument uses finite-dimensional rank
conditions and Kronecker substitution to construct a regular sequence
over $\QQ$ with explicit integer coefficients.

\begin{proposition}\label{prop:specialization}
Let $r\ge1$, and let $w_1,\ldots,w_r$ and $v_1,\ldots,v_r$ be positive integers.  Suppose
there exists a homogeneous regular sequence of degrees $v_1,\ldots,v_r$
in $\mathbb C[x_1,\ldots,x_r]$, with $\wt(x_i)=w_i$.
For $\alpha=(\alpha_1,\ldots,\alpha_r)\in\NN^r$, put
$\wt(\alpha)=\sum_i w_i\alpha_i$. Set
\[
 D=\sum_{j=1}^r v_j-\sum_{i=1}^r w_i,
 \qquad K=\max_i w_i,
 \qquad h_d=\#\{\alpha\in\NN^r:\textstyle\sum_iw_i\alpha_i=d\},
\]
and put
\[
 M=\max_{D<d\le D+K}h_d,\qquad B=M+1,\qquad T=M!+2.
\]
Let $N=\sum_{j=1}^r h_{v_j}$ be the number of coefficients to be specified.
Enumerate all pairs $(j,\alpha)$ with $\sum_iw_i\alpha_i=v_j$ by
distinct indices $e(j,\alpha)=0,1,\ldots,N-1$. Then the following
polynomials form a homogeneous regular sequence over $\QQ$:
\begin{equation}\label{eq:explicit-specialization}
 F_j(x)=\sum_{\wt(\alpha)=v_j}
       T^{\,B^{e(j,\alpha)}}x^\alpha\qquad(1\le j\le r).
\end{equation}
\end{proposition}

\begin{proof}
An existing regular sequence gives an Artinian quotient with Hilbert
series
\[
 H(q)=\frac{\prod_j(1-q^{v_j})}{\prod_i(1-q^{w_i})}.
\]
This polynomial has constant term $1$, and
$H(q^{-1})=q^{-D}H(q)$, so its degree is $D\ge0$.
Consequently every homogeneous component of degree greater than $D$
vanishes.  Also $M\ge1$: for any weight $w_i$, the first multiple of
$w_i$ greater than $D$ is at most $D+K$.

Write $R=\QQ[x_1,\ldots,x_r]$ with these weights, and let $R_d$
be its degree-$d$ part. Introduce independent coefficient indeterminates
$c_{j,\alpha}$ and
generic forms
\[
 f_j=\sum_{\wt(\alpha)=v_j}c_{j,\alpha}x^\alpha.
\]
Put $R_e=0$ for $e<0$ and
write $\boldsymbol c=(c_{j,\alpha})$ for the full coefficient tuple.
Set $\mathcal R_d=R_d\otimes_{\QQ}\QQ[\boldsymbol c]$.
Over the coefficient ring $\QQ[\boldsymbol c]$, consider the matrix
of the multiplication map
\[
 \bigoplus_{j=1}^r \mathcal R_{d-v_j}\longrightarrow \mathcal R_d,
 \qquad (a_1,\ldots,a_r)\longmapsto\sum_j a_j f_j,
\]
in monomial bases.  It has $h_d$ rows, and each matrix entry is either
zero or one coefficient indeterminate.  For $D<d\le D+K$ with $h_d>0$,
evaluation at an existing regular sequence makes this map surjective.
Thus at least one $h_d$ by $h_d$ minor $\Delta_d(\boldsymbol c)$ is a nonzero
polynomial.  Every such minor has degree $h_d\le M$, and its integer
coefficients have absolute value at most $h_d!\le M!$, by the determinant
expansion.

Substitute $c_{j,\alpha}=u^{B^{e(j,\alpha)}}$.  Every exponent of an
individual coefficient in a determinant monomial is at most $M<B$.
Uniqueness of base-$B$ expansions therefore implies that distinct
coefficient monomials become distinct powers of $u$.  Hence
$\widehat\Delta_d(u)$ is nonzero, and every coefficient still has
absolute value at most $M!$.

For completeness, any nonzero integer polynomial
$p(u)=a_su^s+\cdots+a_0$ whose coefficient absolute values are at most
$C_0$ is nonzero at every positive integer $u\ge C_0+2$. When $s\ge1$,
\[
 \left|\sum_{i<s}a_i u^i\right|
 \le C_0\frac{u^s-1}{u-1}<u^s\le |a_s|u^s;
\]
the constant case is immediate. Apply this with $C_0=M!$ and $u=T$.
Thus all the selected minors are nonzero at
\eqref{eq:explicit-specialization}. The ideal
$J=(F_1,\ldots,F_r)\subseteq\QQ[x_1,\ldots,x_r]$ therefore
contains every monomial of weight $d$ for $D<d\le D+K$.

In fact $J$ contains every monomial of weight greater than $D$.
Induct on its weight.  If a monomial has weight $d>D+K$, remove any
variable dividing it.  The remaining monomial has weight at least
$d-K>D$ and smaller than $d$, so it is in $J$ by induction.  Multiplying
back proves the assertion.  Hence $\QQ[x_1,\ldots,x_r]/J$ is Artinian.
The $r$ forms are a system of parameters in an $r$-variable polynomial
ring, and therefore form a regular sequence by the Cohen--Macaulay
argument in Appendix~\ref{app:criterion}.
\end{proof}

\begin{corollary}[A bounded finite search]\label{cor:bounded-search}
Under the hypotheses and notation of Proposition~\ref{prop:specialization},
put
\[
 U=\sum_{D<d\le D+K}h_d.
\]
There is a regular sequence $(F_1,\ldots,F_r)$ over $\QQ$ of the
prescribed degrees such that, for every $j$, the coefficient in $F_j$
of each monomial of weight $v_j$ lies in $\{1,\ldots,U+1\}$.
Thus a finite search of this coefficient grid, checking the ranks of
the multiplication matrices in degrees $D+1,\ldots,D+K$, is guaranteed
to find a regular sequence.
\end{corollary}

\begin{proof}
For each degree with $h_d>0$, choose one of the nonzero minors
$\Delta_d$ from the proof of Proposition~\ref{prop:specialization}.
Their product $P$ is nonzero and has total degree at most $U$, hence
degree at most $U$ in each coefficient variable.  A nonzero polynomial
with that bound cannot vanish on the whole grid
$\{1,\ldots,U+1\}^N$: this follows by induction on the number of variables,
using the fact that a nonzero univariate polynomial of degree at most
$U$ has at most $U$ roots.  At a point where $P$ is nonzero, all required
matrices have full row rank.  The last paragraph of the proposition's
proof then gives regularity.  The search can check matrix ranks directly
and need not identify the chosen minors in advance.
\end{proof}

For the example with weights $4,5,6$ and degrees $10,11,12$, the numbers
$h_{19},\ldots,h_{24}$
are $2,4,3,4,3,5$, so $U=21$ and the bound is $22$.  The supplementary material gives a regular sequence for these
weights and degrees with all allowed coefficients equal to $1$.

\end{document}